\documentclass[10pt,reqno]{amsart}
\usepackage{amsmath}
\usepackage{amscd,amsthm,amsfonts,amsopn,amssymb,mathrsfs}
\usepackage{epsfig,hyperref}

\theoremstyle{remark}

\def\XXint#1#2#3{{\setbox0=\hbox{$#1{#2#3}{\int}$ }
\vcenter{\hbox{$#2#3$ }}\kern-.6\wd0}}

\allowdisplaybreaks

\parskip=\medskipamount

\def\vp{\varphi}

\def\arrowk{^\to{\kern -6pt\topsmash k}}
\def\arrowK{^{^\to}{\kern -9pt\topsmash K}}
\def\arrowr{^\to{\kern-6pt\topsmash r}}
\def\arrowvp{^\to{\kern -8pt\topsmash\vp}}
\def\arrowf{^{^\to}{\kern -8pt f}}
\def\arrowg{^{^\to}{\kern -8pt g}}
\def\arrowu{^{^\to}a{\kern-8pt u}}
\def\arrowt{^{^\to}{\kern -6pt t}}
\def\arrowe{^{^\to}{\kern -6pt e}}
\def\tk{\tilde{\kern 1 pt\topsmash k}}
\def\barm{\bar{\kern-.2pt\bar m}}
\def\barN{\bar{\kern-1pt\bar N}}
\def\barA{\, \bar{\kern-3pt \bar A}}

\def\Re{\text{Re\,}}
\begin{document}
\title{GIBBS MEASURE EVOLUTION IN RADIAL NONLINEAR WAVE AND SCHR\"ODINGER EQUATIONS ON THE BALL\\ \ \\MESURES DE GIBBS ET \'EQUATIONS NON-LIN\'EAIRES DES ONDES ET SCHR\"ODINGER SUR LA BOULE}
\date{\today}

\author{Jean Bourgain}
\address[J.~Bourgain]{School of Mathematics, Institute for Advanced Study, 1 Einstein Drive, Princeton, NJ 08540.}
\email{bourgain@math.ias.edu}
\author{Aynur Bulut}
\address[A.~Bulut]{School of Mathematics, Institute for Advanced Study, 1 Einstein Drive, Princeton, NJ 08540.}
\email{abulut@math.ias.edu}
\thanks{The research of J.B. was partially supported by NSF grants DMS-0808042 and DMS-0835373 and the research of A.B. was supported by NSF under agreement Nos. DMS-0635607 and DMS-0808042.}
\maketitle

\begin{abstract}
We establish new results for the radial nonlinear wave and Schr\"o-\\dinger equations on the ball in $\Bbb R^2$ and $\Bbb R^3$, for random initial data.
More precisely, a well-defined and unique dynamics is obtained on the support of the corresponding Gibbs measure.
This complements results from \cite{B-T1,B-T2} and \cite {T1,T2}.
\end{abstract}

\renewcommand{\abstractname}{R\'esum\'e}
\begin{abstract}
On d\'emontre des r\'esultats nouveaux sur les \'equations des ondes et l'\'equation de Schr\"odinger radiale sur la boule dans $\Bbb R^2$ et $\Bbb R^3$ pour
conditions initiales al\'eatoires. Plus exactement, on 
\'etabl\'e une dynamique bien-d\'efinie et unique sur le support de la measure de Gibbs.  Ceci compl\'emente des r\'esultats de \cite{B-T1,B-T2} et \cite {T1,T2}.
\end{abstract}

\bigskip
\section*{Version fran\c caise abr\'egr\'ee}

On consid\`ere les \'equations non-lineaires (radiale et d\'efocusante) des ondes (NLW) et Schr\"odinger (NLS) sur la boule $B$ dans $\Bbb R^2$ et $\Bbb R^3$
\begin{align*}
(\partial^2_t-\Delta)w +|w|^\alpha w= 0\qquad\qquad&\text{(NLW)}\\
(i\partial_t+\Delta)u -|u|^\alpha u =0 \qquad\qquad &\text{(NLS)}
\end{align*}
ainsi que leures versions tronc\'ees (en introduisant un projecteur $P_N$ sur $[e_1, \ldots, e_N]$ o\`u $\{e_n\}_{n\geq 1}$ sont les fonctions propres de Dirichlet sur $B$) et
les mesures de Gibbs correspondantes.
On \'etablie des estim\'ees espace-temps et une dynamique unique quand $N\to\infty$, dans les mod\`eles (NLW) en dimension 3 pour $\alpha<4$ (le cas $\alpha<3$
\'etant trait\'e dans \cite {B-T1, B-T2}), et (NLS) en dimension 2, $\alpha$ arbitraire (voir \cite {T2} pour le cas $\alpha<4$) et en dimension 3 pour $\alpha=2$.
\bigskip

\section{The equations and the Gibbs measure}

Denote $B=B_d$ the unit ball in $\Bbb R^d$.
We consider the defocusing nonlinear wave (NLW) and nonlinear Schr\"odinger (NLS) equation
\begin{align}
&\left\lbrace\begin{array}{rl}
(\partial^2_t-\Delta)w+|w|^\alpha w&= 0\qquad\qquad\\
(w, \partial_tw)|_{t=0}&=(f_1, f_2)
\end{array}\right.\tag 1\\
&\left\lbrace\begin{array}{rl}
(i\partial_t+\Delta)u-|u|^\alpha u&=0 \qquad\qquad\\
u|_{t=0}&=\phi
\end{array}\right.\tag 2
\end{align}
on the spatial domain $B$ with Dirichlet boundary conditions and with radial initial data.
Thus $(f_1, f_2)$ is real valued and radial in (1), $\phi$ is a radial complex valued function in (2).
It is convenient to rewrite (1) as a first order equation in $t$, introducing the complex function $u=w+i(\sqrt{-\Delta})^{-1} \partial_t w$.
Then (1) turns into the equation
\begin{align}
\begin{cases}
(i\partial_t -\sqrt{-\Delta})u+(\sqrt{-\Delta})^{-1} (|\text{Re}\, u|^\alpha \text{Re}\, u)=0\\
u|_{t=0} =\phi=f_1+i(\sqrt{-\Delta})^{-1} f_2.
\end{cases}
\tag 3
\end{align}
Both (2), (3) are Hamiltonian equations taking the respective forms $i u_t=\frac {\partial H}{\partial\bar u}$ and
$iu_t=(\sqrt {-\Delta})^{-1} \frac {\partial H}{\partial\bar u}$ with Hamiltonians
\begin{align}
H(\phi) =\int_B|\nabla\phi|^2 +\frac 2{2+\alpha} \int_B |\phi|^{\alpha+2}\tag 4
\end{align}
and
\begin{align}
H(\phi)=\int_B|\nabla\phi|^2+\frac 2{2+\alpha} \int_B|\text{Re}\, \phi|^{\alpha+2}.\tag 5
\end{align}
  
Denote $\{e_n\}_{n\geq 1}$ the radial Dirichlet eigenfunctions of $B$ and
$$
P_N\phi =\sum^N_{n=1} \phi_n e_n
$$
the projection operator.
The `truncated' equations
\begin{align}
(i\partial_t +\Delta) u -P_N (|u|^\alpha u)=0\tag 6
\end{align}
\begin{align}
(i\partial_t -\sqrt{-\Delta}) u- P_N\big( (\sqrt{-\Delta})^{-1} (|\Re u|^\alpha \Re u)\big) =0\tag 7
\end{align}
where $\displaystyle u(t)=\sum^N_{n=1} u_n(t)e_n$ are globally wellposed in time and correspond to finite dimensional Hamiltonian models.
The Gibbs measure
\begin{align}
\mu_G^{(N)} (d\phi) = e^{-H(\phi)}\prod^N_1 d^2\phi\tag 8
\end{align}
is invariant under their respective flow.

\section{Statement of the main results}

Our results are the continuation of those obtained in \cite {B-T1,B-T2} and \cite{T1,T2}, as we address various cases that were not treated in these papers.

We consider random initial data given by a Gaussian process
\begin{align}
\phi_\omega =\sum_{n=1}^N \frac {g_n(\omega)}{n\pi} e_n\tag 9
\end{align}
with $\{g_n\}_{n\geq 1}$ independent normalized complex Gaussian random variables.
The free measure $\mu_F^{(N)}$ induced by the map $\omega\mapsto \phi_\omega$ allows to re-express the Gibbs measure as
\begin{align}
\mu_G^{(N)} = e^{-\frac{1}{\alpha+2}\int|\phi|^{\alpha+2}}\mu_F^{(N)}.\tag {10}
\end{align}

Thus, the Gibbs measure is a weighted version of the free measure and has the advantage of being preserved under the flow.
This fact is crucial in the papers cited above and also in the results discussed here.
Note that $\phi_\omega\in H^{\frac 12-}(B)$ almost surely (a.s.).
Fixing $\phi=\phi_\omega$ and considering the truncated solutions $u^N_\phi =u^N$, $u^N|_{t=0} = P_N\phi$ (which are well-defined globally in time), there are
two natural issues.
The first is to establish space-time regularity estimates on $u^N$ that are uniform in $N$.
The second is to prove that for $N\to \infty$, the sequence $\{u^N\}$ converges to a unique limit.
Of course, these properties are only valid a.s. in $\omega$.

\vspace{0.2in}

\noindent {\bf Theorem 2.1} ($3D$ NLW){\bf.}{\it

Let $\alpha<4$.
For almost all $\omega$, the solutions $u^N$ of (7), $u^N|_{t=0}=P_N(\phi_\omega)$ satisfy
\begin{align}
\sup_N\Vert u^N(t)-e^{it\sqrt{-\Delta}} (P_N\phi)\Vert_{H_x^s}<\infty\tag {11}
\end{align}
for all $s<\frac {5-\alpha} 2$ and $t\in\Bbb R$.

Moreover, considering $u^N$ as random variables in $\omega$, the sequence $\{u^N\}$ converges in mean in the space $C_{t<T}H^s_x$ for $s<\frac 12$, $T<\infty$ arbitrary.
}

The case $\alpha<3$ is covered by Theorem 1 in \cite{B-T2}.

\vspace{0.2in}

\noindent {\bf Theorem 2.2} ($2D$ NLS){\bf .}{\it

Let $\alpha \in 2\Bbb Z_+$ be arbitrary and $u^N$ the solutions of (6), $u^N|_{t=0} =P_N(\phi_\omega)$.
Then the sequence $\{u^N\}$ converges in the mean in the space $C_{t<T} H^s_x$ for $s<\frac 12, T<\infty$.
}

The assumption $\alpha\in 2\Bbb Z_+$ is not essential, and more general sufficiently smooth defocusing nonlinearities may be handled as well.  The subquintic case was treated in \cite {T2}.

\vspace{0.2in}

\noindent {\bf Theorem 2.3} ($3D$ NLS){\bf .}{\it 

Let $d=3$ and consider equation (6) with $\alpha =2$.
The solutions $u^N$, $u^N|_{t=0}=P_N(\phi_\omega)$ converge in the mean in the space $C_{t<T} H_x^s$, $s<\frac 12$.
}

\bigskip
\noindent
\section{Comments on the proofs}

As in the many earlier works, the arguments are a combination of probabilistic and harmonic analysis techniques; see for instance the classical works \cite{B0,B1,B3} on this topic.  We only comment on the proof of Theorem $2.3$, which is by far the most delicate.

The starting point is Duhamel's formula on a fixed time interval $[0,T]$
\begin{align}
u^N(t) =u(t) = e^{it\Delta} (P_N\phi)+ i\int_0^t e^{i(t-\tau)\Delta} P_N(u|u|^2)(\tau) d\tau.\tag {12}
\end{align}

The spaces $X_{s,b} =X_{s,b}([0,T])$ are defined in the usual way, see \cite{B4} where they were introduced.  Let $s\geq 0, b\geq 0$.
For functions $f$ on $B_3\times [0, T]$, admitting a representation of the form
\begin{align}
f(x, t)=\sum^\infty_{n=1}\Big[\int^\infty_{-\infty} f_{n, \lambda} \, e^{2\pi i\lambda t} d\lambda\Big] e_n(x) \text { for } x\in B_3\,\, ,\,\, 0\leq t\leq T\tag {13}
\end{align}
where
\begin{align}
\Big(\sum_n\int n^{2s} (1+|n^2-\lambda|^{2b}) |f_{n, \lambda}|^2d\lambda\Big)^{\frac 12}<\infty\tag {14}
\end{align}
we define $\Vert f\Vert_{s, b}$ as the inf (14) over all representations (13).

With these notations, it follows that for $\frac 12< b<1$
\begin{align}
\Big\Vert\int_0^t e^{i(t-\tau)\Delta} f(\tau)d\tau\Big\Vert_{s, b} \leq C\Big(\sum_n \int \frac {n^{2s}|f_{n, \lambda}|^2}
{(1+|n^2-\lambda|^2)^{1-b}}d\lambda\Big)^{\frac 12}. \tag {15}
\end{align}
The inclusions $X_{0, \frac 14+} \subset L_x^{3-} L_t^2$ and $X_{0+, \frac 12+} \subset L^3_xL^4_t$ imply that
$X_{0+, b_1}\subset L_x^3 L_t^{\frac 4{3-4b_1}}$ for $\frac 14< b_1< \frac 12$.
It follows by duality that for $\frac 12 < b<\frac 34$
\begin{align}
(15) \leq\Vert(\sqrt{-\Delta})^{s+} f\Vert_{L_x^{\frac 32} L_t^{\frac 4{5-4b}}}.\tag {16}
\end{align}
From the Gibbs measure conservation under the flow, one derives the a priori inequality (on finite time intervals)
\begin{align}
\Vert(\sqrt{-\Delta})^s u\Vert_{L_x^pL^q_t}<C \ \text { for } \ p<\frac 6{1+2s},\ \,\, q<\infty.\tag {17}
\end{align}
Using (12), (16), (17), it follows that
\begin{align}
\Vert u\Vert_{s, b} < C\ \text { for }\  s<\frac 12,\ \,\, b<\frac 34.\tag {18}
\end{align}

Recall that $u=u_\phi, u|_{t=0} =\phi$ and statements such as (17), (18) require exclusion of small-measure $\phi$-sets.
We do not elaborate on the quantitative aspects of these matters here.

Next, in order to establish convergence properties for $N\to \infty$, let $N\geq N_0$ and estimate using (12) and the preceding
\begin{align}
\nonumber &\Vert u^N -u^{N_0}\Vert_{0, b}\\
&\hspace{0.2in}\leq \Big\Vert\int_0^t e^{i(t-\tau)\Delta} (P_{N_0} u^N|P_{N_0}u^N|^2
-u^{N_0}|u^{N_0}|^2) (\tau)d\tau\Big\Vert_{0, b}+N_0^{-\frac 14}.
\tag {19}
\end{align}

Denoting $u_1=u^{N_0}-P_{N_0} u^N$ and $u_2, u_3$ factors $u_{N_0}, P_{N_0} u^N$, the integrand in (19) leads to trilinear expressions of the form
\begin{align}
\nonumber &\sum_{n, n_1, n_2, n_3} \Big[\int^t_0 \widehat{u_1(\tau)} (n_1) \overline{\widehat{u_2(\tau)}(n_2)} \, \widehat {u_3(\tau)} (n_3) e(-n^2\tau)d\tau\Big]\\
&\hspace{1.8in}\cdot\Big(\int e_n  e_{n_1}  e_{n_2} e_{n_3} dx\Big) e_n e(n^2 t)\tag {20}
\end{align}
where
\begin{align}
\Big|\int e_n e_{n_1} e_{n_2} e_{n_3}\Big|\leq  C\min (n, n_1, n_2, n_3).\tag {21}
\end{align}
Our analysis of (20) is based on arguments closely related to those in \cite{B1}.
We first break up (20) in dyadic regions $n\sim N, n_i\sim N_i(i=1, 2, 3)$ and distinguish the contributions
\begin{align}
|n^2- n_1^2 +n_2^2-n_3^2|\geq \min (N, N_1, N_2+N_3)^{\frac 1{100}}\tag {22}
\end{align}
and
\begin{align}
|n^2- n_1^2 +n_2^2 -n_3^2|<\min (N, N_1, N_2+N_3)^{\frac 1{100}}.\tag {23}
\end{align}

The contribution (22) is handled using $X_{s, b}$-spaces and inequalities of the type
\begin{align}
\iint \bar v u_1 \bar u_2 u_3 dxdt\lesssim\Vert v\Vert_{\sigma_1, 1-b-}\Vert u_1\Vert_{\sigma_2, \frac 12-} \Vert u_2\Vert_{\sigma_3, \frac 34-}
\Vert u_3\Vert_{L^{6-}_x L^ q_t}\tag {24}
\end{align}
with $\sigma_1, \sigma_2\geq 0, \sigma_1+\sigma_2>0, \sigma_3=\frac 12- $ or $\sigma_1=\sigma_2 =0, \sigma_3>\frac 12$ and $q>\frac 4{3-4b}$.

Contributions from (23) are evaluated using further probabilistic considerations, in the spirit of \cite{B1} and exploiting the random nature of $u_2, u_3$.

The most significant terms are
\begin{align}
\sum_n\Big[\int^t_0\widehat {u_1(\tau)}(n)\Big(\sum_m|\widehat{u(\tau)}(m)|^2 \Big(\int e_n^2 e^2_m\Big)\Big)
e(-n^2\tau)\Big] e_n e(n^2t).\tag {25}
\end{align}
Replacing the inner sum $u$ by the free solution $e^{it\Delta}\phi =\sum \frac {g_n(\omega)}n e_n e(n^2t)$ leads to an expression of the form
\begin{align}
\sum_{n<N_0} \log n\Big[\int^t_0 \widehat {u_1(\tau)}(n) e(-n^2\tau)d\tau\Big] e_n e(n^2t).\tag {26}
\end{align}
In order to obtain a contractive estimate in $u_1$, the presence of the $\log n$ factors requires to restrict $t\in [0, T]$, with $T\sim\frac 1{\log N_0}$.

Moreover, the norm $\Vert \cdot \Vert _{0, b}$ has to be slightly weakened to a norm $\vert\vert\vert \cdot \vert\vert\vert_{0, b}$ by allowing in addition to (13), (14) also expressions
$$
f_1(x, t) =\psi(t) \Big[\sum_n b_ne_n e(n^2 t)\Big]
$$
where 
$$
\vert\vert\vert f_1\vert\vert\vert_{0, b} =(\Vert\psi\Vert_\infty +\Vert\psi\Vert_{H^{1/2}})\Big(\sum|b_n|^2\Big)^{\frac 12}<\infty.
$$
Note that the logarithmic divergency above is barely compatible with the error term in (19).

\noindent
{\bf Remark.}
An alternative approach of interest would be to apply the normal forms approach on finite time intervals (cf. \cite{B2}) in order to make reductions of the
Hamiltonian by suitable symplectic transformations.

\end{document}